\documentclass[a4paper,12pt]{article}
\usepackage{amsmath,amsthm,amssymb,latexsym,epic}
\usepackage{graphicx,enumerate}

\newtheorem{theorem}{Theorem}
\newtheorem{lemma}[theorem]{Lemma}
\newtheorem{corollary}[theorem]{Corollary}

\newtheorem{example}[theorem]{Example}
\usepackage[all]{xy}
\newcommand{\IS}{{\mathcal {IS}}}
\newcommand{\dom}{{\mathrm {dom}}}
\newcommand{\ran}{{\mathrm {im}}}
\newcommand{\stran}{{\mathrm {stim}}}
\newcommand{\D}{{\mathcal {D}}}
\newcommand{\J}{{\mathcal {J}}}
\renewcommand{\H}{{\mathcal {H}}}

\begin{document}
\title{On three approaches to conjugacy in semigroups}
\author{Ganna Kudryavtseva and Volodymyr Mazorchuk}
\date{}

\maketitle

\begin{abstract}
We compare three approaches to the notion of 
conjugacy for semigroups, the first one via the transitive
closure of the $uv\sim vu$ relation, the second one via an action
of inverse semigroups on themselves by partial transformations,
and the third one via characters of finite-dimensional representations.
\end{abstract}

\section{Introduction}\label{s1}

Let $G$ be a group. Then it acts on itself via conjugation, 
$x\mapsto gxg^{-1}$, $x,g\in G$; the orbits of this action
define an equivalence relation on $G$ and are called the 
{\em conjugacy classes}. There are alternative equivalent
ways to define conjugate elements. To describe the first one we note
that $y=gxg^{-1}$ can be written as $y=(gx)g^{-1}$ and in this
case $x=g^{-1}(gx)$. Conversely, if $x=uv$ and $y=vu$ for some
$u,v\in G$, then $y=vxv^{-1}$ and thus $x$ and $y$ are conjugate.
Hence the elements $x$ and $y$ are conjugate in $G$ if and only if 
there exist $u,v\in G$ such that $x=uv$ and $y=vu$.

Suppose now that the group $G$ is 
finite and let $\varphi:G\rightarrow \mathbf{GL}(V)$
be a finite-dimensional complex representation of $G$. Then the traces of
the linear operators $\varphi(x)$ and $\varphi(g)\varphi(x)\varphi(g)^{-1}$
coincide for all $x,g\in G$. Conversely, as the characters of irreducible
representations of $G$ form a basis in the space of 
functions constant on conjugacy classes, it follows that $x,y\in G$
are conjugate if and only if for each finite-dimensional representation
$\varphi$ of $G$ the traces of the linear operators $\varphi(x)$ and 
$\varphi(y)$ coincide. 

The aim of the present note, which is inspired by \cite{Hi}, is to extend 
these results to some classes of semigroups. The most traditional 
approach to the notion of conjugacy for semigroups is the one using 
the so-called $G$-conjugacy, defined as follows:  the elements $x,y$ of 
a semigroup $S$, in which $G$ is the group of units, are said to be 
{\em $G$-conjugate}, which is denoted by $x\sim_G y$, provided 
that $x=gyg^{-1}$ for some $g\in G$.
This notion was studied from different points of view by many authors, 
see for example \cite{Pu,Li,KM}. However, this approach is not unique.
Another approach, which will be discussed in detail in the next
section, comes from the equivalence relation generated by the 
(non-transitive in general) relation on $S$, which relates the element
$uv$ to the element $vu$ for all $u,v\in S$. This notion, which has roots
in the study of free monoids, see \cite[11.5]{La},  was also studied in
\cite{GK,KM,KM2,K} for various special classes of semigroups. In the present
paper we show that for some classes of semigroups the latter notion of
conjugacy admits, just like for groups, alternative descriptions via 
the action of a semigroup on itself, and/or via characters of 
finite-dimensional representations.

As usually, we denote {\em Green's relations} on a semigroup by 
$\mathcal{J}$, $\mathcal{D}$, $\mathcal{R}$, $\mathcal{L}$ and 
$\mathcal{H}$. If $\varphi:S\rightarrow \mathrm{End}_{\mathbb{C}}(V)$
is a representation of a semigroup $S$, then the {\em character} of 
$\varphi$ is the function $\chi_{\varphi}:S\rightarrow \mathbb{C}$ for 
which $\chi_{\varphi}(s)$ is defined as the trace of $\varphi(s)$, $s\in S$.
A semigroup $S$ is called {\em group-bound} provided that for each $x\in S$
there exists $k\in\{1,2,\dots\}$ such that $x^k$ lies in a subgroup of
$S$. If $S$ is group-bound, then for $x\in S$ we denote by $e_x$ the
idempotent of the subgroup containing $x^k$ as above. 
For a partial transformation $a$ we denote by $\dom(a)$ and $\ran(a)$ 
the {\em domain} or the {\em image}  of $a$ respectively.
\vspace{0.5cm}

\noindent
{\bf Acknowledgment.} This work was done during the visit of the
first author to Uppsala University, which was supported by the 
Royal Swedish Academy of Sciences  (KVA) and the Swedish Foundation for
International Cooperation in Research and Higher Education (STINT). For the
second author the research was partially supported by the Swedish Research
Council. The financial support of the Swedish Research Council, KVA and 
STINT and the hospitality of Uppsala University are gratefully acknowledged.

\section{Three approaches to the notion of conjugacy 
for semigroup}\label{s2}

\subsection{The first approach: via the transitive closure of the
relation $uv\sim vu$}\label{s2.1}

Let $A$ be a non-empty alphabet and $A^*$ the corresponding free monoid.
Two elements $x, y\in A^*$ are called {\em conjugate} if 
there are $u,v\in A^*$ such that $x=uv$ and $y=vu$, or, equivalently, 
if there is $u\in A^*$ such that $ux=yu$ (see e.g. \cite[11.5]{La}). 
These two relations, however, do not coincide on non-free semigroups
in general (for example if $S$ contains a zero element, then the second
relation obviously degenerates). Hence if $S$ is a semigroup, we will call
$x,y\in S$ {\em primarily conjugate} and denote this fact by $x \sim_p y$ 
if there are $u,v\in S^1$ such that $x=uv$ and $y=vu$. The relation 
$\sim_p$ is reflexive and symmetric, while not transitive in 
general, in contrast to what one has for groups and free semigroups. 
Let $\sim$ be the transitive
closure of $\sim_p$. We call two elements $x,y\in S$ {\em conjugate} 
if $x\sim y$. It is easily verified that for monoids $x\sim_G y$ 
always implies $x\sim y$. The converse inclusion is not true
in the general case, see for example \cite{GK}. The relation $\sim$ was
studied for many classes of semigroup, see for example \cite{La,GK,KM,KM2,K}.

\subsection{The second approach: via an action of a semigroup 
on itself}\label{s2.2}

Let $S$ be  an inverse group-bound semigroup, and $\geq$ denote the 
natural partial order on $S^1$, that is $a\geq b$ if 
and only if there is an idempotent $e\in S$ such 
that $b=ae$ (see also \cite[5.2]{how} for alternative definitions). 
If $e,f\in S$ are idempotents then $e\geq f$ is equivalent to the 
equality $ef=fe=f$. Let $a\in S^{1}$, $x\in S$. Set
\begin{equation}\label{eq_action}
a\cdot x =\left\lbrace\begin{array}{l} axa^{-1}, {\text{ if }} a^{-1}a\geq e_x;\\
{\text{not defined, otherwise}.}
\end{array}\right.
\end{equation}

It is obvious that $x\mapsto a\cdot x$, $x\in S$, is a partial one-to-one transformation of $S$. Moreover, the following lemma shows that 
\eqref{eq_action} in fact defines an action of $S^1$ on $S$.

\begin{lemma}\label{lem:well_def}
Let $a,b\in S^{1}$, $x \in S$. Then $ba\cdot x$ is defined if and only if $a\cdot x$ and $b\cdot(a\cdot x)$ are defined. Moreover, in the case when $ba\cdot x$ is defined we have the equality $ba\cdot x=b\cdot(a\cdot x)$.
\end{lemma}

\begin{proof}
Using the Preston-Wagner representation (see e.g. \cite[Theorem~1.20]{CP}) 
we can assume that the semigroup $S$ is a subsemigroup of some inverse 
symmetric semigroup  $\IS(X)$. For a group-bound element $a\in\IS(X)$ 
the element $e_a$ is an idempotent 
acting identically on the set consisting of all $t\in X$ which belong 
to cycles of $a$ (see \cite[5.1]{GM} for the case of finite $X$, the case
of infinite $X$ can be treated easily using similar arguments). This 
set is called {\em the stable image} of $a$ and is denoted by 
$\stran(a)$ (see for example \cite{GM}). 
The condition $a^{-1}a\geq e_x$ is then equivalent to 
$\dom(a)\supseteq \stran(x)$.

Observe that whenever the condition $\dom(a)\supseteq\stran(x)$ holds 
there is a bijection (given by $(t_1,\dots, t_l)\mapsto 
(a(t_1),\dots, a(t_l))$) from cycles of $x$ to cycles of $axa^{-1}$. 
This, in particular, implies that $a(\stran(x))=\stran(axa^{-1})$.

Suppose that $a,b,x\in \IS(X)$ are group-bound elements and  
that $a\cdot x$ and $b\cdot(a\cdot x)$ are defined. Since $a\cdot x$ 
is defined we have $\dom(a)\supseteq\stran(x)$, which implies that 
\begin{equation}\label{eq:aux1}
\ran(a)\supseteq a(\stran(x))=\stran(axa^{-1}).
\end{equation}
In addition, since $b\cdot(a\cdot x)$ is defined we have the inclusion 
\begin{equation}\label{eq:aux2}
\dom(b)\supseteq\stran(axa^{-1}).
\end{equation}
From~\eqref{eq:aux1} and~\eqref{eq:aux2} it follows that
$\dom(b)\cap\ran(a)\supseteq a(\stran(x))$, whence
\begin{displaymath}
\dom(ba)=a^{-1}(\dom(b)\cap\ran(a))\supseteq a^{-1}(a(\stran(x)))=\stran(x),
\end{displaymath}
which means that $ba\cdot x$ is defined.

Suppose now that $ba\cdot x$ is defined. This means that $\dom(ba)\supseteq\stran(x)$ which implies that 
$a^{-1}(\ran(a)\cap\dom(b))\supseteq \stran(x)$. This, in turn, gives us 
\begin{equation}\label{eq:aux3}
\dom(b)\cap\ran(a)\supseteq a(\stran(x)).
\end{equation}
Therefore, $\dom(a)\supseteq \stran(x)$,  and thus $a(\stran(x))=\stran(axa^{-1})$. This and~\eqref{eq:aux3} imply that $\dom(b)\supseteq\stran(axa^{-1})$, and hence both $a\cdot x$ 
and $b\cdot(a\cdot x)$ are defined.

That $ba\cdot x=b\cdot(a\cdot x)$  whenever $ba\cdot x$, $a\cdot x$ and
$b\cdot(a\cdot x)$ are defined follows from the definition of our action 
and the fact that for an inverse semigroup the operation of taking the 
inverse of a given element is an anti-involution.
\end{proof}

Let $S$ be an inverse semigroup and  $x,y\in S$. Set $x\approx_p y$ if 
there is some $a\in S^1$ such that $y=a\cdot x$ or $x=a\cdot y$. This relation 
is reflexive because $x=1\cdot x$. It is also symmetric by definition.  
Let $\approx$ be the transitive closure of $\approx_p$. We will call this 
relation {\em the relation of conjugacy in the action sense}.

\subsection{The third approach: via characters of finite-dimensional
representations}\label{s2.3}

As we have mentioned in the Introduction, two elements $x,y$ of a finite
group $G$ are conjugate if and only if for every finite-dimensional 
complex representation $\varphi$ of $G$ the equality
$\chi_{\varphi}(x)=\chi_{\varphi}(y)$ holds. Let now $S$ be a 
semigroup and $x,y\in S$.  We will call $x$ and $y$ {\em conjugate in 
the character sense} and  denote this fact by $x\equiv y$ provided that 
for every finite-dimensional complex representation $\varphi$ of $S$ 
we have the equality $\chi_{\varphi}(x)=\chi_{\varphi}(y)$. 

\section{Comparing the three approaches}\label{s3}

\begin{theorem}\label{th1}
Let $S$ be a regular group-bound semigroup with finite $\D$-classes. 
Then for $x,y\in S$ we have  $x\sim y$ if and only if $x\equiv y$.
\end{theorem}

\begin{proof}
If $a,b\in S$ are such that $a\sim_p b$, then $a=uv$ and $b=vu$ 
for some $u,v\in S$. Let $\varphi$ be a finite dimensional complex representation of $S$. Then the traces of the linear 
operators $\varphi(u)\varphi(v)$ and $\varphi(v)\varphi(u)$ are equal
(by a standard exercise in linear algebra), and thus the necessity of 
the claim follows.

Suppose now that $a,b\in S$ and $a\equiv b$. Since $S$ is regular and 
group-bound, from \cite[Corollary~6]{K} for all $x,y\in S$ it follows 
that $x\sim y$ if and only if  $xe_x\sim ye_y$. Having this in mind, 
it is enough  to prove that $ae_a\sim be_b$.  

Let us first show that 
$a\equiv ae_a$ and $b\equiv be_b$. It is of course enough to prove the 
first formula. 
Let $\varphi:S\to  \mathrm{End}_{\mathbb{C}}(V)$ be a finite-dimensional 
representation of $S$. For $t\in \mathbb{C}$ set
$V_t=\{v\in V\,:\, (\varphi(a)-t)^{\dim(V)} v=0\}$.
From the definition it follows that the linear operator $\varphi(e_a)$
it the projection of $V$ onto $U=\oplus_{t\neq 0}V_t$ with the kernel $V_0$.
In particular, the actions of $\varphi(a)\varphi(e_a)$ and $\varphi(a)$ on
$U$ coincide, and both $\varphi(a)\varphi(e_a)$ and
$\varphi(a)$ act nilpotently on $V_0$. This implies that
$\chi_{\varphi}(a)=\chi_{\varphi}(ae_a)$ and hence $a\equiv ae_a$.
As a consequence we also obtain $ae_a\equiv be_b$.

Let us now show that the elements $ae_a$ and $be_b$ belong to the same 
$\D$-class of $S$. To do this we have to recall the 
construction of some induced modules for semigroups, which follows
closely \cite[5.4]{CP}. 

For $x\in S$ denote by $D_x$, $L_x$, and 
$H_x$ the  $\D$-class, $\mathcal{L}$-class or $\H$-class of $x$, respectively.
Let $e\in S$ be an idempotent and $\varphi: H_e\rightarrow
\mathbf{GL}(W)$ be a finite-dimensional representation of the maximal
subgroup $H_e$. Let $H_1=H_e$, $H_2,\dots,H_k$ be the list of all
$\H$-classes in $L_e$. Fix some $a_i$ in each $H_i$ such that
$a_i$ is an idempotent if $H_i$ is a subgroup. 

Let $s\in S$ and $i\in\{1,2,\dots,k\}$. Assume that $sa_i\in D_e$ and
let $a'_i$ be some inverse to $a_i$ (which exists as $S$ is regular). 
Then $sa_i=(sa_ia'_i)a_i$, where $sa_ia'_i$ must belong to $D_e$
as well. As $D_e$ is finite, from \cite[Lemma~1.3.3]{how} it follows 
that $sa_i\in L_e$. This implies that we have exactly two possibilities:
either $sa_i\not\in L_e$ (and hence $sa_i\not\in D_e$) or 
$sa_i=a_js'$ for some uniquely determined $j\in\{1,2,\dots,k\}$ and 
$s'\in H_e$. For $i\in\{1,2,\dots,k\}$ let $W^{(i)}$ denote a copy 
of $W$. Then we
can consider the vector space $\overline{W}=\oplus_{i=1}^k W^{(i)}$
and for every $s\in S$ define a linear operator on $\overline{W}$
as follows: for $v\in W^{(i)}$ set
\begin{equation}\label{eqeq}
sv=
\begin{cases}
0, & sa_i\not\in L_e;\\
s'v\in W^{(j)},&  sa_i=a_js'\text{ as above}. 
\end{cases}
\end{equation}
This defines a representation $\overline{\varphi}:S\rightarrow
\mathrm{End}_{\mathbb{C}}(\overline{W})$. Observe that, by 
\cite[Theorem~2.17]{CP}, for $s\in H_e$ we have $sa_i\in L_e$ if and only 
if $a_i$ is an idempotent. Moreover, in this case for $a_j$ and
$s'$ from \eqref{eqeq} we have $a_j=a_1=e$ 
and $s'=s$ by \cite[Lemma~2.14]{CP}. This implies that for $s\in H_e$ we
have $\chi_{\overline{\varphi}}(s)=m\chi_{\varphi}(s)$, where $m>0$ is 
the number of idempotents in $L_e$.

Suppose that $ae_a$ and $be_b$ do not belong to the same $\D$-class of $S$.
Then they do not belong to the same $\J$-class either for on 
group bound semigroups Green's relations $\D$ and $\J$ coincide,
see \cite[Theorem~1.2.20]{Hi1}. Hence, without loss of generality we can 
assume that $ae_a\not\in S be_b S$.  The element $ae_a\in H_{e_a}$ is
a group element. Let $\varphi: H_{e_a}\rightarrow \mathbf{GL}(W)$
be some finite-dimensional complex representation such that 
$\chi_{\varphi}(ae_a)\neq 0$. Such representation exists since characters
of irreducible representations form a basis in the space of class functions.
Then $\chi_{\overline{\varphi}}(ae_a)=m\chi_{\varphi}(ae_a)\neq 0$, while
$\overline{\varphi}(be_b)=0$ as $ae_a\not\in Sbe_b S$. This contradicts
$ae_a\equiv be_b$ and proves that the elements $ae_a$ and $be_b$ do
belong to the same  $\D$-class of $S$.

Since $ae_a\D be_b$ and both are group elements, from \cite[Theorem~2.20]{CP}
it follows that there exist a pair, $t$ and $t'$, of mutually inverse 
elements in the 
same $\D$-class such that $tbe_bt'\H ae_a$. Let $\psi$ be a 
finite-dimensional representation of $S$. Then we have
\begin{equation}\label{eq:a}
\chi_{\psi}(tbe_bt^{-1})=\chi_{\psi}(be_bt^{-1}t)=
\chi_{\psi}(be_b)=\chi_{\psi}(ae_a),
\end{equation}
where the second equality follows from $be_bt^{-1}t=be_b$, which, in turn,
follows from \cite[Lemma~2.14 and Theorem~2.17]{CP}. In particular,
the characters of $tbe_bt^{-1}$ and $ae_a$ coincide in all cases when
$\psi=\overline{\varphi}$, where $\varphi$ is an irreducible 
representation of $H_{ae_a}$. As on elements from $H_{ae_a}$ 
the character of $\overline{\varphi}$ differs from that of $\varphi$ only 
by a  non-zero constant (the number of idempotents in $L_{ae_a}$, 
see above), it follows that the characters of $tbe_bt^{-1}$ and 
$ae_a$ coincide for every irreducible representation of $H_{ae_a}$. 
Therefore $tbe_bt^{-1}$ and $ae_a$
are conjugate as elements of $H_{ae_a}$, in particular
$tbe_bt^{-1}\sim ae_a$. As $be_b=(be_bt^{-1})t$, the
elements $be_b$ and $tbe_bt^{-1}$ are primarily conjugate. This implies
that $ae_a\sim be_b$ and completes the proof.
\end{proof}

\begin{theorem}\label{th2}
Let $S$ be a group-bound inverse semigroup and $x,y\in S$. Then the following conditions are equivalent:
\begin{enumerate}[(a)]
\item \label{i1} $x\sim y$;
\item \label{i2} $x\approx y$;
\item \label{i3} there is $z\in S$ such that 
$z=a\cdot x=b\cdot y$ for some $a,b\in S^1$.
\end{enumerate}
\end{theorem}

\begin{proof}
The implication \eqref{i3}$\Rightarrow$\eqref{i2} is obvious.

Let us prove the implication \eqref{i1}$\Rightarrow$\eqref{i3}. 
Using the Preston-Wagner representation we again think of $S$ as of 
a group-bound subsemigroup of the inverse symmetric semigroup $\IS(X)$. 
Suppose that $x,y\in S$ are such that $x\sim y$. We take $z=xe_x$. 
First we observe that $\dom(e_x)=\stran(x)$, and thus $e_x\cdot x$ is 
defined. Now we observe that $e_x\cdot x=e_xxe_x=xe_x$. Therefore 
$x\approx xe_x$ and  we can take  $a=e_x$. Analogously one shows that
$e_y\cdot y=ye_y$.

Since  $x\sim y$ it follows from \cite[Corollary~6]{K} 
that $xe_x\sim ye_y$ and 
$xe_x\D ye_y$. By \cite[Theorem~2.20]{CP} there exists 
$t\in L_{e_y}\cap R_{e_x}$ such 
that $xe_x=tye_yt^{-1}$. Since $\dom(t)=\dom(e_y)=\stran(ye_y)$, it follows 
that $t\cdot ye_y$ is defined and equals $xe_x$. Hence 
$te_y\cdot y$ is defined and  equals $xe_x$ by Lemma~\ref{lem:well_def}. 
We take $b=te_y$ and the implication \eqref{i1}$\Rightarrow$\eqref{i3} follows.

Finally, we show that \eqref{i2}$\Rightarrow$\eqref{i1}. 
Again, using the Preston-Wagner representation we think of $S$ as of 
a group-bound subsemigroup of the inverse symmetric semigroup $\IS(X)$. 
Let $x,y\in S$ be such that $y=a\cdot x=axa^{-1}$ for some $a\in S$
such that $a^{-1}a\geq e_x$. Then $\dom(a)\supseteq \stran(x)$ and
it follows that $a$ induces a bijection between $\stran(x)$ and
$\stran(y)$. Let $\hat{a}=ae_x\in S$ denote the restriction of $a$ to
$\stran(x)$. Then we have $\hat{a}^{-1}\hat{a}=e_x$,
$\hat{a}\hat{a}^{-1}=e_y$ by definition. Moreover, it also
follows that $ye_y=\hat{a}xe_x\hat{a}^{-1}$ and
$xe_x=\hat{a}^{-1}ye_y\hat{a}$. Hence
\begin{displaymath}
xe_x=\hat{a}^{-1}ye_y\hat{a}\sim
ye_y\hat{a}\hat{a}^{-1}=ye_y.
\end{displaymath}
Applying \cite[Corollary~6]{K} we obtain $x\sim y$, which completes the proof.
\end{proof}

From Theorems \ref{th1} and \ref{th2} we immediately obtain the following corollaries.

\begin{corollary}
Let $S$ be a finite regular semigroup and $a,b\in S$. Then $a\sim b$ 
if and only if $a\equiv b$.
\end{corollary}

\begin{corollary}
Let $S$ be a finite inverse semigroup and $a,b\in S$. Then $a\sim b$ 
if and only if $a\approx b$ if and only if $a\equiv b$.
\end{corollary}

It would be interesting to extend the definition of the relation
$\approx$ to other classes of semigroups. We finish the paper with a
an example, which shows that Theorem~\ref{th1} is not true for 
non-regular semigroups.

\begin{example}\label{ex1}
{\rm
Let $S=\langle x\,:\, x^2=x^3\rangle$. Then the relation $\sim$ on
$S$ is trivial (i.e. all conjugacy classes contain exactly one element).
At the same time if $\varphi:S\rightarrow \mathrm{End}_{\mathbb{C}}(V)$
is a finite-dimensional complex representation of $S$, then 
$\varphi(x^2)$ is a projection (as $x^2$ is an idempotent), so 
the same arguments as in the proof of Theorem~\ref{th1} show that
$\chi_{\varphi}(x)=\chi_{\varphi}(x^2)$, implying $x\equiv x^2$.
}
\end{example}

\vspace{1cm}

\noindent
V.M.: Department of Mathematics, Uppsala University, SE 471 06,
Uppsala, SWEDEN, e-mail: {\em mazor\symbol{64}math.uu.se}
\vspace{0.5cm}

\noindent
G.K.: Department of Mechanics and Mathematics, Kyiv Taras Shev\-chen\-ko
University, 64, Volodymyrska st., 01033, Kyiv, UKRAINE,
e-mail: {\em akudr\symbol{64}univ.kiev.ua}

\end{document}